\documentclass[12pt,a4paper]{article}
\usepackage[cp1251]{inputenc}
\usepackage[T2A]{fontenc}
\usepackage[english]{babel}

\usepackage{amsmath,amsfonts,amssymb,mathrsfs,amsopn,indentfirst,amscd}

\usepackage{graphicx}
\usepackage{amsthm}
\usepackage{hyperref}

\usepackage{amssymb}
\usepackage{amsmath}

\newtheorem*{corollary*}{Corollary}

\begin{document}

\title{On a Simplex Inscribed in a Ball}

\author{Mikhail Nevskii\footnote{ Department of Mathematics,  P.\,G.~Demidov Yaroslavl State University,\newline Sovetskaya str., 14, Yaroslavl, 150003, Russia, \newline
               mnevsk55@yandex.ru,  orcid.org/0000-0002-6392-7618 }
}               
      
\date{May 21, 2026}
\maketitle

\begin{abstract}

\smallskip
 Let $B_n$ be the $n$-dimensional unit ball given by the inequality 
$\|x\|\leq 1$, where $\|x\|$ is the standard Euclid norm in  ${\mathbb R}^n$. 
For an $n$-dimensional nondegenerate simplex 
$S$, we denote by $E$ the ellipsoid of minimum volume which contains $S$. Suppose $S\subset B_n$,
$0\leq m\leq n-1$.
Let $G$ be any $m$-dimensional face of $S$ and let $H$ be the opposite $(n-m-1)$-dimensional face. Denote by $g$ and $h$  the centers of gravity of $G$ and $H$ respectively. Define   $y$ as the intersection point of the line passing from $g$ to $h$ with the boundary of $E$. Let us call the face $G$ suitable if $y\in B_n.$ 
Earlier it was  proved that each simplex $S\subset B_n$ has a~suitable face of any dimension $\leq n-1$.
We~show the following. Let $S$ be inscribed in $B_n$. If some vertex of $S$ is suitable, then there exists a suitable face of any dimension $\leq n-1$ which contains this vertex.

\smallskip 

Keywords: $n$-dimensional ball, $n$-dimensional simplex, minimum volume ellipsoid, linear interpolation

\smallskip
MSC: 41A05, 52B55, 52C07

\end{abstract}

\section{Introduction}\label{nev_s1}

Hereafter $n\in{\mathbb N}$. For
$x=(\xi_1,\ldots,\xi_n), y=(\eta_1,\ldots,\eta_n)\in {\mathbb R}^n$,
by ${\langle x,y\rangle}$ we denote the~standard scalar product in
${\mathbb R}^n$:
$$\langle x,y\rangle=\xi_1\eta_1+\ldots+\xi_n\eta_n.$$
Let $\|x\|=\sqrt{\langle x,x\rangle}$ be the Euclid norm in ${\mathbb R}^n$.  We denote by $B_n$
the unit $n$-dimensional ball given by the inequality  
$\|x\|\leq 1.$

F. John \cite{nevskii_john} proved that every $n$-dimensional convex body 
contains a single ellipsoid of maximum volume, and also gave a characterization of convex bodies for which such an~ellipsoid is the unit ball
$B_n$ (for more details see, e.g., \cite{nevskii_ball_90}, \cite{nevskii_ball_97}). From John's theorem,  the analogous statement follows, which characterizes the only ellipsoid of minimum volume containing a given convex body.

We will consider an ellipsoid of minimum volume containing a given nondegenerate simplex. We will call such an ellipsoid {\it the minimal ellipsoid of a simplex}. Clearly, the~minimal ellipsoid of a simplex is circumscribed around this simplex.
The center \linebreak of the ellipsoid coincides with the center of gravity of the simplex.
The minimal ellipsoid of a simplex is~a~Euclidean ball if and only if this simplex is regular. This is equivalent to the well-known property that of all simplices contained in a ball the maximum volume has the~regular simplex inscribed in this ball 
(see, e.\,g., \cite{nevskii_fejes_tot},
\cite{nevskii_slepian},
\cite{nevskii_vandev}).

In this paper we will consider some questions related to the mutual location \linebreak of~an~$n$-dimensional ball
and the minimal ellipsoid of a simplex contained in this ball.

Consider a nondegenerate $n$-dimensional simplex $S$ with vertices $x_1,\ldots, x_{n+1}$. We~denote by $E$ the minimum volume ellipsoid
containing $S$. Let $m$ be an integer, $1\leq m\leq n$, and let $J\subset \{1, \ldots, n+1\}$ be an arbitrary set consisting of~$m$ elements.
Let us put in~correspondence to $J$  the point $y_J\in E$ defined as follows.  Assume
$g_J$ is
 the center of gravity of the $(m-1)$-dimensional face of $S$ containing vertices $x_j, j\in J$, and $h_J$ is the center of~gravity of the $(n-m)$-dimensional face containing the remaining $n+1- m$ vertices. Then $y_J$ is the intersection point of the line passing through $g_J$ and $h_J$ in direction from the first point to the second with~the boundary of the ellipsoid $E$.

\medskip

 {\bf Theorem 1.}
  {\it Suppose $S\subset B_n$. Then for any $m\in\{1,\dots,n\}$  there exists a set
$J\subset \{1, \ldots, n+1\}$ consisting of $m$  numbers such that $y_J\in B_n$.
}

\medskip
If for a set $J$ of $m$ numbers the point $y_J$ belongs to $B_n$, then the $(m-1)$-dimensional face of $S$ containing the vertices $x_j, j\in J$,
will be called {\it suitable}. Thus, Theorem 1 states that every simplex 
$S\subset B_n$ has a suitable face of any  dimension $\leq n-1.$

The origin of this subject is related to polynomial interpolation of functions of many variables. Namely, Theorem 1  was proved  in \cite{nevskii_mz_23} as an important element in solving the problem of finding the minimum norm of a projector in linear interpolation
 on $B_n$.
The application of this theorem allows us to establish  that the minimum  $C(B)$ operator norm corresponds to a  projector with  nodes at the vertices of a regular simplex inscribed in the ball.
The exact value
$\theta_n(B_n)$ of the minimum norm of a  projector is 
$$\theta_n(B_n)=\max\{\psi(a_n),\psi(a_n+1)\},$$
 where 
$$\psi(t):=\frac{2\sqrt{n}}{n+1}\Bigl(t(n+1-t)\Bigr)^{1/2}+
\left|1-\frac{2t}{n+1}\right|, \quad 0\leq t\leq n+1,$$
$$a_n:=\left\lfloor\frac{n+1}{2}-\frac{\sqrt{n+1}}{2}\right\rfloor.$$
We have
$$\sqrt{n}\leq \theta_n(B_n)\leq \sqrt{n+1}.$$
Moreover,  $\theta_n(B_n)$ $=$ $\sqrt{n}$ only for $n=1$, and
 $\theta_n(B_n)=\sqrt{n+1}$
 if and only if
$\sqrt{n+1}$ is an integer. For application
of Theorem~1, the important case is $m=k_n$, where $k_n$ coincides with that of the numbers $a_n$ and $a_n+1$ on~which $\psi(t)$ takes a larger value. For more details, see ~\cite{nevskii_mz_23}, \cite{nevskii_uh_19}, \cite{nevskii_21}.

At the same time, this geometric result seems to be interesting in itself.  In the present paper, in a different way than in \cite{nevskii_mz_23}, we prove 
the following proposition.

\medskip
 {\bf Theorem 2.} {\it Let any simplex $S$ be inscribed in $B_n$.
  If some vertex of $S$ is suitable, then for each
   $m\in\{1,\ldots, n-1\}$
the simplex $S$ has a suitable  $m$-dimensional face containing this vertex.
}

\medskip
The proof of this theorem is given in Section \ref{nev_s3}.

\section{Preliminaries}\label{nev_s2}

Let us make some preliminary remarks.
The center of gravity of  $S$ and the center of  minimal ellipsoid $E$ coincide and are located at the point 
$$c=\frac{1}{n+1}\sum_{j=1}^{n+1} x_j.$$
 Denote by $r$ the ratio of the distance from the center of gravity of the regular simplex to~the center of gravity of the $(m-1)$-dimensional face
to the radius of the circumscribed sphere.
 The number $r$ is easily found by considering a regular simplex inscribed in $B_n$. 
 The center of gravity of such a simplex is $c=0$, so 
$$
\langle x_i,\sum_{j=1}^{n+1} x_j\rangle =\langle x_i,nc\rangle=\langle x_i,0\rangle=0,$$
wherefrom
$$\langle x_i,x_j\rangle=-\frac{1}{n}, \quad i\ne j.$$ 
Consequently, we have
$$r^2=\left\| \frac{1}{m}\sum_{j=1}^m x_j\right\|^2=\frac{1}{m^2}\left(\sum_{j=1}^m \|x_j\|^2 +2\sum_{1\leq i<j\leq m}
\langle x_i,x_j\rangle\right)=$$
$$=\frac{1}{m^2}\left (m-2\cdot\frac{m(m-1)}{2n}\right)=\frac{1}{m^2}\cdot\frac{mn-m^2+m}{n}
=\frac{n-m+1}{mn}.$$
Thus
$$r=
\sqrt{\frac{n-m+1}{mn}}.$$  
Put
\begin{equation}\label{varrho}
\varrho=\frac{1}{r}= \sqrt{\frac{mn}{n-m+1}}.
\end{equation}
By $r^\prime$ and $\varrho^\prime$, we denote the analogs of  numbers $r$ and $\varrho$ for an $m$-dimensional face: 
\begin{equation}\label{r_prime_varrho_prime}
r^\prime=\sqrt{\frac{n-m}{(m+1)n}}, \quad \varrho^\prime=\frac{1}{r^\prime}= \sqrt{\frac{(m+1)n}{n-m}}.
\end{equation}

Consider an arbitrary set $J$ of $m$ indices
$j.$ The center of gravity of an $(m-1)$-dimensional face with~vertices $x_j,$  $j\in J,$ is the point
$$g_J=\frac{1}{m}\sum_{j\in J} x_j.$$
Since an ellipsoid is affinely equivalented to a ball, the point $y_J$ has the form
$$y_J=c+\varrho\left(c-g_J\right)=\left(1+\varrho\right)c-\varrho g_J.$$
Therefore, 
\begin{equation}\label{y_J}
\|y_J\|^2=\langle y_J, y_J \rangle =(1+\varrho)^2\|c\|^2-2\varrho(1+\varrho)\langle c,g_J\rangle +
\varrho^2\|g_J\||^2.
\end{equation}
For any $m$, the sum $Q$ of coefficients in front of the squares $\|x_i\|^2$ in the right-hand side of 
\eqref{y_J} is equal to $1$:
$$Q= (1+\varrho)^2\cdot \frac{1}{(n+1)^2}\cdot (n+1)-2\varrho(1+\varrho)\cdot\frac{1}{m(n+1)}\cdot m+
\varrho^2\cdot\frac{1}{m^2}\cdot m=$$
$$=\frac{(1+\varrho)^2-2\varrho(1+\varrho)}{n+1}+
\frac{\varrho^2}{m}
=\frac{1}{n+1}+\varrho^2\left(\frac{1}{m}-\frac{1}{n+1}\right)=$$
$$=\frac{1}{n+1}+\frac{mn}{n-m+1}\cdot\frac{n+1-m}{m(n+1)}=\frac{1}{n+1}+\frac{n}{n+1}=1.$$
We have applied \eqref{varrho}.

Also note that 
a suitable vertex from the condition of Theorem 2 exists for any simplex $S\subset B_n$. This follows from Theorem 1 but can also be obtained directly.  The center of~gravity of the simplex $c$ lies on the hyperplane
$\langle c,c-x\rangle=0,$ so some vertex $x_j$ belongs to the half-space $\langle c,c-x\rangle\leq 0$. This vertex is just  a suitable one. Indeed, for the set $J=\{j\}$ containing a single
index, we have $\varrho=1, g_J=x_j, y_J=2c-x_j$. Hence,
$$\\|y_J\|^2=\langle y_J, y_J \rangle=4\|c\|^2-4\langle c,x_j\rangle +\|x_j\||^2=4\langle c,c-x_j\rangle+\|x_j\|^2\leq \|x_j\||^2\leq 1.
$$
Thus, Theorem 2 is some generalization of Theorem 1 for an inscribed $S$.

\section{Proof of Theorem 2}\label{nev_s3}

By the condition, the simplex $S$
with vertices $x_i$ is inscribed
in $B_n,$ i.\,e., $\|x_i\|=1.$ 

The case $n=1$ is trivial since then $S=E=B_1$. Suppose $n>1.$
Let the vertex $x_1$ of~the~simplex is suitable, i.\,e., for $J=\{1\}$ the point $y_J$ belongs to $B_n$.
Since $y_J=2c-x_1$ we have
\begin{equation}\label{norm_2c_minus_x_1_quadr}
\|y_J\|^2=4\|c\|^2-4\langle c,x_1\rangle +\|x_1\|^2\leq 1.
\end{equation}
As noted above, the sum of the terms with $\|x_i\|^2$ included in \eqref{norm_2c_minus_x_1_quadr} is equal to $1$. Therefore,
 the inequality $\|y_J\|^2\leq 1$
means that in the middle part of \eqref{norm_2c_minus_x_1_quadr} the sum of scalar products  $\langle x_i,x_j\rangle$ for $i< j$
 with~the corresponding coefficients  is nonpositive.
 It is easy to see that \eqref{norm_2c_minus_x_1_quadr} reduces to the inequality
$$\frac{8}{(n+1)^2}\sum\limits_
{ 2\leq i<j\leq n+1}\langle x_i,x_j\rangle +\left(\frac{8}{(n+1)^2}- \frac{4}{n+1}\right)\sum_{j=2}^{n+1} \langle x_1,x_j\rangle\leq 0, $$.
which is equivalent to
 \begin{equation}\label{condition_for_x_1}
\sum\limits_
{ 2\leq i<j\leq n+1}\langle x_i,x_j\rangle \leq \frac{n-1}{2}\sum_{j=2}^{n+1} \langle x_1,x_j\rangle.
\end{equation}

 Fix $m$, $1\leq m \leq n-1$.
Let us prove that \eqref{condition_for_x_1} implies the existence of pairwise different indices $2\leq j_1,\ldots,j_m\leq n+1$ such that  
 for $J^\prime=\{1,j_1,\ldots,j_m\}$ the point $y_{J^\prime}$ belongs to~$B_n.$ This will mean
 that the $m$-dimensional face of $S$ containing  vertices $x_1, x_{j_1},\ldots, x_{j_m}$ is  suitable. The number of all sets of the considered type is $n\choose{m}$. Therefore, suffice~it to~show that
 \begin{equation}\label{sum_sqnorms_y_J_prime_leq_n_m}
 \sum_{J^\prime} \|y_{J^\prime}\||^2 \leq {n\choose m}.
 \end{equation}

 Since $|{J^\prime}|=m+1$, we have
$$\|y_{J^\prime}\|^2 =(1+\varrho^\prime)^2\|c\|^2-2\varrho^\prime(1+\varrho^\prime)\langle c,g_{J^\prime}\rangle +
{\varrho^\prime}^2\|g_{J^\prime}\|^2,
$$
where $\varrho^\prime$ is defined by the equality \eqref{r_prime_varrho_prime}. Summing over all  
$J^\prime=\{1,j_1,\ldots,j_m\}$,
we can write down
\begin{equation}\label{sum_y_J_prime}
\sum_{J^\prime} \|y_{J^\prime}\|^2={n\choose m}(1+\varrho^\prime)^2 \|c\|^2-2\varrho^\prime(1+\varrho^\prime) 
\langle c, \sum_{J^\prime} g_{J^\prime}\rangle +{\varrho^\prime}^2 \sum_{J^\prime} \|g_{J^\prime}\||^2.
\end{equation}
Notice that
$$(m+1) \sum_{J^\prime} g_{J^\prime}={n\choose m} x_1+{{n-1}\choose {m-1}}\sum_{j=2}^{n+1} x_j=
\left[ {n\choose m}-{{n-1}\choose{m-1}}\right]x_1+{{n-1}\choose {m-1}}(n+1)c,$$
whence
$$\sum_{J^\prime} g_{J^\prime}
=\frac{ {n\choose m}-{{{n-1}\choose{m-1}}} }{m+1} x_1+{{n-1}\choose {m-1}}\frac{n+1}{m+1}c.
$$
Let us continue \eqref{sum_y_J_prime} by applying the latter equality:
$$\sum_{J^\prime} \|y_{J^\prime}\|^2
=\|c\|^2\left[{n\choose m}(1+\varrho^\prime)^2-2\varrho^\prime(1+\varrho^\prime){{n-1}\choose {m-1}}\frac{n+1}{m+1}\right]-
$$
\begin{equation}\label{long_eq_for_sum_sqnorms_y_J_prime}
-2\varrho^\prime(1+\varrho^\prime) 
\frac{ {n\choose m}-{{n-1}\choose{m-1}} }{m+1}  \langle c,x_1\rangle 
+{\varrho^\prime}^2  \sum_{J^\prime} \|g_{J^\prime}\|^2.
\end{equation}
From the above, the contribution to this sum of terms with squared norms $\|x_i\|^2$
 is equal to $n\choose m$, i.\,e., the number of $\|y_{J^\prime}\|^2$. Therefore, \eqref{sum_sqnorms_y_J_prime_leq_n_m}
 is equivalent to the fact that the sum of scalar products $\langle x_i,x_j \rangle$ for $i<j$
with corresponding coefficients entering $\sum \|y_{J^\prime}\||^2$ is~nonpositive.
Let us distinguish two groups of scalar products: products $\langle
x_i,x_j\rangle$ for $2\leq i<j\leq n+1$ and products $\langle x_1,x_j\rangle$ for $2\leq j\leq n+1$, 
and then  rewrite  inequality~\eqref{sum_sqnorms_y_J_prime_leq_n_m} in the form 
 \begin{equation}\label{sums_x_i_x_j_x_1_x_j_ineq}
a\sum\limits_
{  2\leq i<j\leq n+1}\langle x_i,x_j\rangle \leq b\sum_{j=2}^{n+1} \langle x_1,x_j\rangle.
\end{equation}
The coefficients $a$ and $b$ are determined from the right-hand side of \eqref{long_eq_for_sum_sqnorms_y_J_prime},
their transformation is rather labor-consuming.  Here we give only basic formulae. From \eqref{long_eq_for_sum_sqnorms_y_J_prime} we have:
$$a=\frac{2}{(n+1)^2}
\left[{n\choose m}(1+\varrho^\prime)^2-2\varrho^\prime(1+\varrho^\prime){{n-1}\choose {m-1}}\frac{n+1}{m+1}\right]+\frac{2}{(m+1)^2}{{n-2}\choose{m-2}}{\varrho^\prime}^2$$
(excepting the last summand when $n=2$),
$$b=-\frac{2}{(n+1)^2}
\left[{n\choose m}(1+\varrho^\prime)^2-2\varrho^\prime(1+\varrho^\prime){{n-1}\choose {m-1}}\frac{n+1}{m+1}\right]+$$
$$+2\varrho^\prime(1+\varrho^\prime)
\frac{ {n\choose m}-{{n-1}\choose{m-1}} }{(m+1)(n+1)}-
\frac{2}{(m+1)^2}{{n-1}\choose{m-1}}{\varrho^\prime}^2.$$
Making use of the equalities
$${{n-1}\choose m}=\frac{n-m}{n}{n\choose m}, \quad
{{n-1}\choose {m-1}}=\frac{m}{n}{n\choose m}, \quad
{{n-2}\choose {m-2}}=\frac{m(m-1)}{n(n-1)}{n\choose m},$$
as well as $${\varrho^\prime}^2= \frac{(m+1)n}{n-m}$$
(see \eqref{r_prime_varrho_prime}), we can represent $a$ and $b$ in the following form:
\begin{equation}\label{a_formulae}
a=\frac{4}{(n+1)^2(m+1)} {n\choose m} \left[ \frac{n-m}{n} \varrho^\prime+\frac{2m^2+n^2-3mn+m-n}{(n-m)(n-1)}\right],
\end{equation}
\begin{equation}\label{b_formulae}
b=\frac{2}{(n+1)^2(m+1)} {n\choose m} \left[ \frac{(n-m)(n-1)}{n} \varrho^\prime+\frac{2m^2+n^2-3mn+m-n}
{n-m}\right].
\end{equation}
It follows that
\begin{equation}\label{b_a_eq_n-1_2}
\frac{b}{a}=\frac{n-1}{2}.
\end{equation}

Thus, the coefficients in  inequalities \eqref{condition_for_x_1} and \eqref{sums_x_i_x_j_x_1_x_j_ineq}
turn out to be proportional. Since the vertex $x_1$ is suitable, 
we have \eqref{condition_for_x_1} 
and with it also \eqref{sums_x_i_x_j_x_1_x_j_ineq}. The latter, as already noted, is equivalent to
\eqref{sum_sqnorms_y_J_prime_leq_n_m}. So, the average value of $\||y_{J^\prime}\|^2$ over
all sets $J^\prime=\{1,j_1,\ldots,j_m\}$ is~not greater than 1. Consequently, there exists a set $J^\prime$
of the considered form for which $ \|y_{J^\prime}\|^2\leq 1.$ This means that some $m$-dimensional face
containing $x_1$ is suitable.

Finally, we remark that  equality \eqref{b_a_eq_n-1_2} can be obtained in a simpler way,
without deriving the explicit  formulae \eqref{a_formulae}--\eqref{b_formulae} for  coefficients
 $a$ and $b$. Clearly, 
 $a$ and $b$ do not depend on a simplex inscribed into the ball. Let us choose
$S$ as a regular simplex inscribed in $B_n$. 
As we have noted, $\langle x_i,x_j\rangle=-1/n, \ i\ne j.$ 
 Since in this case the minimal ellipsoid $E$ coincides with $B_n$, all $\|y_{J^\prime}\|$ are equal to $1$ and
inequality~\eqref{sum_sqnorms_y_J_prime_leq_n_m} becomes an equality. Hence, the~equivalent relation
\eqref{sums_x_i_x_j_x_1_x_j_ineq}  is also an equality:  
$$a\cdot\frac{(n-1)n}{2}\cdot\left(-\frac{1}{n}\right)=b \cdot \left(-\frac{1}{n}\right)\cdot n.$$
This  gives  \eqref{b_a_eq_n-1_2}.
The proof is complete.
\hfill$\Box$

\section{Concluding Remarks}\label{nev_s4}

The proposition generalizing Theorem 2 can be formulated as follows.

\medskip
{\it Let $n \geq 2$ and $0\leq m\leq n-2$. If the $m$-dimensional face $G$
of a simplex $S$ inscribed in $B_n$ is~ suitable, then $S$ has a suitable face of any dimension
$d\in\{m+1,\ldots,n-1\}$
containing~$G$.}

\medskip
Theorem 2 corresponds to the case $m=0$ when $G$ is a vertex. The other extreme case
$m=n-2$ is trivial because any $(n-1)$-dimensional face of a simplex $S\subset B_n$ is suitable.
Indeed, for such a face, the set $J$ consists of $n$ numbers, $\varrho=n$ (see \eqref{varrho})
 and $y_J$ coincides with the vertex
$x_k$ where $k\not\in J:$
$$y_J=\left(1+\varrho\right)c-\varrho g_J=\sum_{j=1}^{n+1} x_j-\sum_{j\in J} x_j=x_k.$$
However, the validity of this proposition for $0<m<n-2$ remains unclear.

\end{document}